\documentclass[12pt]{article}
\usepackage{amsmath, amssymb}

\newtheorem{Theorem}{Theorem}[section]

\newtheorem{Remark}[Theorem]{Remark}

\newtheorem{Example}[Theorem]{Example}

\def\R{\mathbb{R}}
\def\E{\mathbb{E}}
\def\N{\mathbb{N}}

\numberwithin{equation}{section}

\begin{document}

\noindent {\Large \bf A Burgers-KPZ Type Parabolic Equation
\par\noindent
for the Path-Independence of the Density
\par\noindent
of the Girsanov Transformation}
\bigskip
\bigskip
\par\noindent
{\bf Aubrey Truman$^{b}$, Feng-Yu Wang$^{a,b}$, Jiang-Lun Wu$^{b}$
and Wei Yang$^{b}$}
\smallskip
\par\noindent
{\footnotesize{${}^{a}$School of Mathematical Sciences, Beijing
Normal University, Beijing 100875, China}}
\par\noindent
{\footnotesize{${}^b{}$Department of Mathematics, Swansea
University, Singleton Park, Swansea, SA2 8PP,}}
\par\noindent
 {\footnotesize{\,\,  and Wales Institute of Mathematical and Computational Sciences, UK}}
\par\noindent
{\footnotesize{E-mail: a.truman@swansea.ac.uk,
f.y.wang@swansea.ac.uk,

\qquad j.l.wu@swansea.ac.uk, mawy@swansea.ac.uk}}

\bigskip
\begin{abstract}
Let $X_t$ solve the multidimensional It\^o's stochastic differential
equations on $\R^d$
$$dX_t=b(t,X_t)dt+\sigma(t,X_t)dB_t$$
where $b:[0,\infty)\times\R^d\to\R^d$ is smooth in its two
arguments, $\sigma:[0,\infty)\times\R^d\to\R^d\otimes\R^d$ is smooth
with $\sigma(t,x)$ being invertible for all
$(t,x)\in[0,\infty)\times\R^d$, $B_t$ is $d$-dimensional Brownian
motion. It is shown that, associated to a Girsanov transformation,
the stochastic process
$$\int^t_0\langle(\sigma^{-1}b)(s,X_s),dB_t\rangle+\frac{1}{2}\int^t_0|\sigma^{-1}b|^2(s,X_s)ds$$
is a function of the arguments $t$ and $X_t$ (i.e.,
path-independent) if and only if $b=\sigma\sigma^\ast\nabla v$ for
some scalar function $v:[0,\infty)\times\R^d\to\R$ satisfying the
time-reversed KPZ type equation
$$\frac{\partial}{\partial t}v(t,x)=-\frac{1}{2}\left[\left(Tr(\sigma\sigma^\ast\nabla^2v)\right)(t,x)
+|\sigma^\ast\nabla v|^2(t,x)\right].$$ The assertion also holds on
a connected complete differential manifold.
\end{abstract}

{\medskip\par\noindent {\bf Mathematics Subject Classification
(2000)}: 60H10, 58J65, 35Q53.
\smallskip\par\noindent
{\bf Key Words and Phrases}: } Stochastic differential equations;
the Girsanov transformation; nonlinear parabolic partial
differential equations of Burgers-KPZ type; diffusion processes and
nonlinear PDEs on differential manifolds.

\section{Introduction and motivations}
The present paper is mainly concerned with a link of It\^o's
stochastic differential equations (SDEs) to nonlinear parabolic
partial differential equations (PDEs) of Burgers-KPZ type, by virtue
of Girsanov transformation. Our result presents a characterization
of the path-independence property for the density process of
Girsanov transformation to SDEs. Our second interest is to establish
such a connection between SDEs and nonlinear PDEs on complete
differential manifolds.

Since the pioneering work of J.M. Burgers in 1930s (cf. e.g.
\cite{Burgers}), Burgers equation -- the simplest nonlinear PDE
$$\frac{\partial}{\partial t}u(t,x)+\lambda
u(t,x)\frac{\partial}{\partial
x}u(t,x)=\nu\frac{\partial^2}{\partial x}u(t,x), \quad
(t,x)\in[0,\infty)\times\R$$ has received a great attention both in
mathematics and physics. Wherein the parameter $\lambda\in\R$
measures the strength of the nonlinearity, $\nu>0$ stands for the
viscosity and the (linear) viscous dissipation term on the right
hand side of the equation is for the sake of softening shock wave
phenomena.

Fix $d\in\mathbb{N}$, let $\mathbb{R}^d$ be the $d$-dimensional
Euclidean space with the inner product being denoted by
$\langle\cdot,\cdot\rangle$. The multidimensional analogue to the
above Burgers equation is the so called higher dimensional Burgers
equation for a vorticity-free velocity field ${\bf
u}:[0,\infty)\times\R^d\to\R^d$ (cf. e.g., \cite{AlbeverioMS}) which
reads the following
$$\frac{\partial}{\partial t}{\bf u}+\lambda({\bf u}\cdot\nabla){\bf u}=\nu\Delta{\bf
u}$$ where $\nabla$ stands for the space gradient, the dot product
${\bf u}\cdot\nabla:=\langle{\bf u},\nabla\rangle$, and
$\Delta:=\nabla\cdot\nabla$ the Laplace operator on $\R^d$.
Nowadays, Burgers equation is significant in the mathematical
modeling of the large scale structure of the universe with
complexity. The equation appears in many fields like aerodynamics,
fluid dynamics (in particular, hydrodynamics), polymers and
disordered systems, turbulence and propagation of chaos, as well as
in shock wave and conservation laws -- to name just a few. Among
many interesting and important investigations, a breakthrough study
has been made by three physicists M.P. Kardar, G. Parisi, and Y.-C.
Zhang (\cite{KPZ}) for modeling the time evolution of the profile of
a growing interface with the name of Kardar-Parisi-Zhang equation,
or in short, KPZ equation. The KPZ equation describes the
macroscopic properties of a wide variety of growth processes, such
as growth by ballistic deposition and the Eden model (cf.
\cite{KSpohn}). For a more mathematical account of the KPZ equation,
the reader is referred to \cite{Handa}. The link of the KPZ equation
to multidimensional Burgers equation can be explicated as follows.
It is a natural assumption that the field ${\bf u}$ is often
generated by a potential function (i.e., the profile)
$u:[0,\infty)\times\R^d\to\R$
$${\bf u}(t,\cdot)=-\nabla u(t,\cdot),\quad t\in[0,\infty)$$
which, from the multidimensional Burgers equation, gives to the
following KPZ equation for $u$
$$\frac{\partial}{\partial t}u(t,x)=\nu\Delta u(t,x)+\frac{\lambda}{2}|\nabla
u(t,x)|^2\, .$$ Clearly, the above KPZ equation describes the
large-distance, long-time dynamics of the growth process specified
by a single-valued height $u(t,x)$ on a substrate $x\in\R^d$. It
reflects the competition between the surface tension smoothing
forces $\nu\Delta u(t,x)$ and $\frac{\lambda}{2}|\nabla u(t,x)|^2$
(the nonlinear term of $u$ represents the tendency for growth to
occur preferentially in the local normal direction to the surface).

When the diffusion coefficient $\sigma\equiv\sigma_0$, a constant, 
very interesting and new links of (stochastic) multidimensional Burgers' 
equations to (stochastic) Hamilton-Jacobi-Bellman (in short, HJB) 
equations and the continuity equation have been thoroughly investigated 
by Truman and Zhao in \cite{TrumanZhao,TrumanZhao1,TrumanZhao2} 
(see also the early works \cite{ElworthyTruman,Truman} for a bridge between 
the diffusion equations and the Schr\"odinger equation, now called the 
Elworthy-Truman formula). In this content, Hamilton-Jacobi continuity 
equations provide the key to obtaining asymptotic expansions in ascending 
powers of $\sigma_0$ for solutions of the corresponding heat (and Schr\"odinger) 
wave functions in this setting. Actually, the iterated Hamilton-Jacobi-continuity 
equations derived there inspired our consideration carried out in the present work.

Nowadays, because of their ubiquity, Burgers equation, the KPZ
equation, and the HJB equations (as well as any of their advances studies) 
maintain a very hot research topic on both theoretical and applied aspects 
in various fields involving disordered systems and non-equilibrium dynamics.
The applied aspect links to many diverse areas ranging from physics,
biochemistry, and climate and ocean studies (cf. e.g.
\cite{Spohn,Woycz}), to economical and financial studies (cf.
\cite{HC1993,HL2004,Stein,CoxLeland,YZ}). There are many works in the
literature devoted to analytic aspect of the equations themselves as
well as to computational aspect (cf. e.g.
\cite{Dafermos,Freidlin,Majda,MT,Smoller,CIL} and references therein).

On the other hand, the theory of SDEs has been very well developed
since the seminal work of the great Japanese mathematician Kiyosi
It\^o in the mid 1940s, cf. \cite{Ito}. Since then, SDEs have
profound impacts on differential geometry and PDEs (cf.
\cite{Elworthy,Elworthy1,IkWa,Malliavin} and most recently \cite{Stroock}).

In recent years, due to the necessity of introducing stochastic
volatility as the measurement of uncertainty in modeling of
financial markets, stochastic differential equations receive a huge
attention from both theoretical and practical aspects, cf. e.g.
\cite{GS,IkWa,Kunita,MaThalm,Oksendal,StrV}. The primary point here
is to model the price dynamics or the wealth growth by utilizing
SDEs, after established a so-called real world probability space
(cf. e.g. the seminal paper \cite{BlackScholes} by Black and
Scholes). It is a pivotal problem to characterize the
path-independence property for certain utility functions in an
equilibrium market. We will give some concrete exposition of this
point after the presentation of our first main result, which is the
motivation of our study from economics and finance.

The object of the present paper is to explore a novel link from
It\^o's SDEs to nonlinear parabolic PDEs of Burgers-KPZ type with
our particular attention to derive such a connection from SDEs to
nonlinear PDEs on differential manifolds. Our results give a
characterization of path-independence of the density of the Girsanov
transformation for SDEs in terms of a nonlinear parabolic PDE of
Burgers-KPZ type.

The rest of the paper is organized as follows. In the next section,
we first give a brief account of the Girsanov transformation for
multidimensional SDEs on $\R^d$, then we formulate our result on the
characterization of path-independence of the Girsanov density and we
give some further account of our result to relevant studies on
path-independence features in economics and finance. We then present
our proof to the theorem and we end Section 2 with an exposition of
the one dimension case. Section 3, the final section, is devoted to
the extension of the connection of SDEs to nonlinear equation of
Burgers-KPZ type on connected complete differential manifolds.

\section{The characterization theorem on $\R^d$}

\subsection{Preliminaries on SDEs and the Girsanov transformation}
Let us start with the general framework of stochastic differential
equations by following \cite{IkWa}. Given a complete probability
space $(\Omega,\mathcal{F},P)$ with a usual filtration
$\{\mathcal{F}_t\}_{t\in[0,\infty)}$. Let
$C([0,\infty),\mathbb{R}^d)$ be the space of all
$\mathbb{R}^d$-valued, continuous functions defined on $[0,\infty)$.
It is known that $C([0,\infty),\mathbb{R}^d)$ is a complete,
separable metric space under the metric
$$\rho(u,v):=\sum^{\infty}_{k=1}2^{-k}(\mbox{max}_{0\le t\le k}|u(t)-v(t)|\wedge1),
\quad u,v\in C([0,\infty),\mathbb{R}^d)\, .$$ We endow the
topological $\sigma$-algebra
$\mathcal{B}(C([0,\infty),\mathbb{R}^d))$ on
$C([0,\infty),\mathbb{R}^d)$ so that
$(C([0,\infty),\mathbb{R}^d),\mathcal{B}(C([0,\infty),\mathbb{R}^d)))$
forms a measurable space, and further we denote by
$\mathcal{B}_t(C([0,\infty),\mathbb{R}^d))$ the sub-$\sigma$-algebra
of $\mathcal{B}(C([0,\infty),\mathbb{R}^d))$ generated by the family
$\{C([0,\infty),\R^d)\ni u\mapsto u(s): 0\le s\le t\}$ for
$t\in[0,\infty)$. As usual, $\mathbb{R}^d\otimes\mathbb{R}^m$ (with
$m\in\mathbb{N}$) stands for the totality of real $d\times m$
matrices (realised alternatively by identifying
$\mathbb{R}^d\otimes\mathbb{R}^m$ with $dm$-dimensional Euclidean
space) endowed with the Hilbert-Schmidt norm
$$|a|=\sqrt{\sum^d_{j=1}\sum^m_{k=1}|a^j_k|^2},
\quad a=(a^j_k)_{d\times m}\in \mathbb{R}^d\otimes\mathbb{R}^m$$ and
$\mathcal{B}(\mathbb{R}^d\otimes\mathbb{R}^m)$ denotes the
topological $\sigma$-algebra on $\mathbb{R}^d\otimes\mathbb{R}^m$.
Moreover, we use the notation $\mathcal{A}^{d,m}$ to denote the
collection of all
$\mathcal{B}([0,\infty))\times\mathcal{B}(C([0,\infty),
\mathbb{R}^d))/\mathcal{B}(\mathbb{R}^d\otimes\mathbb{R}^m)$-measurable
mappings
$$a:[0,\infty)\times
C([0,\infty),\mathbb{R}^d)\to\mathbb{R}^d\otimes\mathbb{R}^m$$ such
that for each $t\in[0.\infty)$, the mapping
$$u\in C([0,\infty),\mathbb{R}^d)\mapsto a(t,u)\in\mathbb{R}^d\otimes\mathbb{R}^m$$
is
$\mathcal{B}_t(C([0,\infty),\mathbb{R}^d))/\mathcal{B}(\mathbb{R}^d\otimes\mathbb{R}^m)$-measurable.

Given $b\in\mathcal{A}^{d,1}$ and $\sigma\in\mathcal{A}^{d,d}$, we
consider the following stochastic differential equation of the
Markovian type for a $d$-dimensional continuous process
$X=(X_t)_{t\in[0,\infty)}$

\begin{equation} \label{SDE}
dX_t=b(t,X_t)dt+\sigma(t,X_t)dB_t, \quad t\ge0
\end{equation}
where
\begin{gather*}
X_t=
\begin{bmatrix} X^1_t\\X^2_t\\ \vdots\\ X^d_t\\ \end{bmatrix},
\quad b=
\begin{bmatrix} b^1\\ b^2\\ \vdots\\ b^d\\ \end{bmatrix},
\quad \sigma=
\begin{bmatrix}
    \sigma^1_1&\cdots & \sigma^1_d\\
    \sigma^2_1& \cdots & \sigma^2_d\\
    \vdots& &  \vdots\\
    \sigma^d_1&\cdots & \sigma^d_d\\
\end{bmatrix},
\quad B_t=
\begin{bmatrix}B^1_t\\B^2_t\\ \vdots\\ B^d_t\\ \end{bmatrix},
\end{gather*}
so equation (\ref{SDE}) in terms of its components is
\begin{equation}\label{SDEa}
dX^j_t=b^j(t,X_t)dt+\sum^d_{k=1}\sigma^j_k(t,X_t)dB^k_t,\quad
j=1,2,...,d
\end{equation}
where $\sigma^j_k$ stands for the $(j,k)$-entry of the $d\times
d$-matrix $\sigma$, for $j,k=1,2,...,d$, and
$B_t=(B^1_t,B^2_t,...,B^d_t)^\ast$ is an $d$-dimensional
$\{\mathcal{F}_t\}_{t\in[0,\infty)}$-Brownian motion. It is well
known, from e.g. \cite{IkWa} (cf. Theorem IV.3.1), that under the
usual conditions of linear growth and locally Lipschitz, to be
precise, the coefficients $b:[0,\infty)\times\R^d\to\R^d$ and
$\sigma:[0,\infty)\times\R^d\to\R^d\otimes\R^d$ satisfy linear
growth and locally Lipschitz condition, $C^1$ with respect to the
first variable, and $C^2$ with respect to the second variable, there
exists a unique solution to equation (\ref{SDE}) with given initial
data $X_0$. By Stroock-Varadhan's martingale problem \cite{StrV},
$X_t$ is associated with the following second order elliptic
differential operator (called the Markov generator)
$$L_tf(x)=\frac{1}{2}\sum^d_{i,j=1}a^{ij}(t,x)
\frac{\partial^2f(x)}{\partial x_i\partial
x_j}+\sum^d_{j=1}b^i(t,x)\frac{\partial f(x)}{\partial x_i}, \quad
f\in C^2(\R^d)$$ with $a(t,x):=\sigma(t,x)\sigma^\ast(t,x)$, where
$\sigma^\ast(t,x)$ stands for the transposed matrix of
$\sigma(t,x)$. In component form,
$a^{ij}(t,x):=\sum^d_{k=1}\sigma^i_k(t,x)\sigma^j_k(t,x)$.

The celebrated Girsanov theorem provides a very powerful
probabilistic tool to solve equation (\ref{SDE}) under the name of
the {\it Girsanov transformation or the transformation of the
drift}. Let $\gamma\in\mathcal{A}^{d,1}$ satisfy the following
condition
$$\E\left[\exp\left(\frac{1}{2}\int_0^t|\gamma(s,X_s)|^2ds\right)\right]
<\infty, \quad \forall t>0.$$ Then, by Girsanov theorem (cf e.g.
Theorem IV 4.1 of \cite{IkWa}),
$$\exp\left(\int^t_0\gamma(s,X_s)dB_s-\frac{1}{2}\int^t_0|
\gamma(s,X_s)|^2ds\right),\quad t\in[0,\infty)$$ is an
$\{\mathcal{F}_t\}$-martingale. Furthermore, for $t\ge0$, we define
$$Q_t:=\exp\left(\int^t_0\gamma(s,X_s)dB_s-\frac{1}{2}
\int^t_0|\gamma(s,X_s)|^2ds\right)\cdot P$$ or equivalently in terms
of the Radon-Nikodym derivative
$$\frac{dQ_t}{dP}=\exp\left(\int^t_0\gamma(s,X_s)dB_s
-\frac{1}{2}\int^t_0|\gamma(s,X_s)|^2ds\right).$$ Then, for any
$T>0$,
$$\tilde{B}_t:=B_t-\int^t_0\gamma(s,X_s)ds,\quad 0\le t\le T$$
is an $\{\mathcal{F}_t\}$-Brownian motion under the probability
$Q_T$. Moreover, $X_t$ satisfies
$$dX_t=[b(t,X_t)+\sigma(t,X_t)\gamma(t,X_t)]dt
+\sigma(t,X_t)d\tilde{B}_t, \quad t\geq 0.$$

One can then discuss comprehensively the {\it existence and
uniqueness as well as the structure of solutions} to the initial
value problem for equation (\ref{SDE}) by appealing the above
argument with suitable choice of $\gamma$. Here we want to explore
such transformation to another link to partial differential
equations.

\subsection{The characterization theorem and its link to economics and finance studies}
From now on in the paper, we assume the coefficient $\sigma$
satisfies that the matrix $\sigma(t,x)$ is invertible, for any
$(t,x)\in[0,\infty)\times\R^d$ (and consequently so is the symmetric
matrix $a(t,x)=\sigma(t,x)\sigma^{\ast}(t,x)$). Moreover, we specify
the above function $\gamma$ by
$$\gamma(t,x)=-(\sigma(t,x))^{-1}b(t,x)$$
so that $b(t,X_t)+\sigma(t,X_t)\gamma(t,X_t)=0$, and hence we
further require $b$ and $\sigma$ satisfy
$$\E\left[\exp\left(\frac{1}{2}\int_0^t|(\sigma(s,X_s))^{-1}
b(s,X_s)|^2ds\right) \right]<\infty, \quad \forall t>0.$$ Thus the
associated probability measure $Q_t$ is determined by
\begin{eqnarray*}
\frac{dQ_t}{dP}&=&\exp\left(-\int^t_0\langle(\sigma(s,X_s))^{-1}
b(s,X_s),dB_s\rangle\right.\\
&& \quad \left.-\frac{1}{2}\int^t_0\big|(\sigma(s,X_s))^{-1}
b(s,X_s)\big|^2ds\right)\,.\end{eqnarray*}

Set
$$\hat{Z}_t:=-\ln\frac{dQ_t}{dP}$$
that is
$$\hat{Z}_t=\int^t_0\langle(\sigma(s,X_s))^{-1}b(s,X_s),
dB_s\rangle+\frac{1}{2}\int^t_0\big|(\sigma(s,X_s))^{-1}
b(s,X_s)\big|^2ds\,.
$$
Clearly, $\hat{Z}_t$ is a one dimensional stochastic process with
the stochastic differential form
$$d\hat{Z}_t
=\frac{1}{2}\big|(\sigma(t,X_t))^{-1} b(t,X_t)\big|^2dt
+\langle(\sigma(t,X_t))^{-1}b(t,X_t),dB_t\rangle\, .$$

We are now ready to state the first main result of this paper. It
gives a necessary and sufficient condition, and hence a
characterization of path-independence of the density $\hat{Z}_t$ of
the Girsanov transformation for SDEs in terms of a nonlinear
parabolic PDE of Burgers-KPZ type. Namely we establish a bridge from
SDE (\ref{SDE}) to a nonlinear parabolic PDE of Burgers-KPZ type.

\begin{Theorem} \label{TH1}
Let $v:[0,\infty)\times\R^d\to\R$ be a scalar function which is
$C^1$ with respect to the first variable and $C^2$ with respect to
the second variable. Then

\begin{eqnarray}\label{RND0}
v(t,X_t) &=& v(0,X_0)+\frac{1}{2}\int^t_0\big|(\sigma(s,X_s))^{-1}
b(s,X_s)\big|^2ds\nonumber \\
&& \quad +\int^t_0\langle(\sigma(s,X_s))^{-1} b(s,X_s),dB_s\rangle
\end{eqnarray}
equivalently,
$$\frac{dQ_t}{dP}=\exp\{v(0,X_0)-v(t,X_t)\},\quad t\in[0,\infty)$$
holds if and only if
\begin{equation}\label{eqv}
b(t,x)=(\sigma\sigma^\ast\nabla v)(t,x),\quad
(t,x)\in[0,\infty)\times\R^d\end{equation} and $v$ satisfies the
following time-reversed KPZ type equation
\begin{equation}\label{KPZ-type}
\frac{\partial}{\partial
t}v(t,x)=-\frac{1}{2}\left\{\big[Tr(\sigma\sigma^\ast\nabla^2v)\big](t,x)
+|\sigma^\ast\nabla v|^2(t,x)\right\}
\end{equation}
where $\nabla^2v$ stands for the Hessian matrix of $v$ with respect
to the second variable.
\end{Theorem}

\begin{Remark}\label{rm0} 
The derived time-reversed KPZ type equation (\ref{KPZ-type}) is contained 
as a special version of the 
stochastic HJB equation derived in \cite{TrumanZhao2}. It is an interesting 
question to see if one can recover a fuller picture of the mathematical physics 
of the stochastic HJB equations in \cite{TrumanZhao,TrumanZhao1,TrumanZhao2} 
by exploiting the argument developed in this paper. We will consider 
this problem in our future work. 
\end{Remark}

Before presenting the proof to Theorem \ref{TH1}, we would like to
give some links of our result to economic and financial studies.

\begin{Remark}\label{rm1}
Recall that in economics and finance studies, a conventional kind of
equilibrium financial market can be characterized by the utility
function of a representative agent (see e.g.,
\cite{CoxLeland,D19981,D19982, DR2003}). Given the probability
measure $P$ as an objective probability in the market model, one can
interpret our process $X_t$ as the wealth (or the assets price) of
the representative agent in a multi-assets market. If the class of
utility functions is one of differentiable, increasing, and strictly
concave time-separable von Neumann-Morgenstern utility functions,
then the representative agent maximizes his/her expected utility
$U$. Cox and Leland in \cite{CoxLeland} show that the
path-independence property is necessary for expected utility
maximization. By path-independence, they mean that the value of a
portfolio will depend only on the assets price at that time, not on
the path followed by the assets in reaching that price (vector).
Namely, the utility function $U$ depends on the state price $X_t$
and time $t$, for $t\geq 0$, that is, the function $U$ is of the
form $U(X_t, t)$. On the other hand, Dybvig and Ross in
\cite{DR2003} show that, in an equilibrium market, the marginal
utility in each state is proportional to a consistent state-price
density function. In addition, in an equilibrium market, there must
exist a risk neutral probability measure $Q$ which is absolutely
continuous with respect to $P$. The Radon-Nikodym derivative
$Z=\frac{dQ}{dP}$ gives the state-price density \cite{HC1993}.
Combining the above $U(X_t,t)$, therefore, the Radon-Nikodym
derivative is also in the form of
$$Z(X_t,t)=\frac{dQ_t}{dP}.$$
Clearly, our Theorem \ref{TH1} presents a necessary and sufficient
condition for the above Radon-Nikodym derivative, hence a
characterization for the path-independence property of the utility
function.
\end{Remark}

\subsection{Proof of Theorem \ref{TH1}} \label{Proof}
We start with the {\it necessity}. Namely, assume that there exists
a scalar function $v:[0,\infty)\times\R^d\to\R$ which is $C^1$ with
respect to the first variable and $C^2$ with respect to the second
variable such that (\ref{RND0}) holds. Then by (\ref{RND0}), we have

\begin{equation} \label{v_t2}
    dv(t,X_t)=\frac{1}{2}\big|(\sigma(t,X_t))^{-1}
    b(t,X_t)\big|^2dt+\langle(\sigma(t,X_t))^{-1} b(t,X_t),dB_t\rangle \, .
\end{equation}

Now by viewing $v(t,X_t)$ as the composition of the deterministic
$C^{1,2}$-function $v:[0,\infty)\times\R^d\to\R$ with the continuous
semi-martingale $X_t$, we can apply It\^o's formula to $v(t,X_t)$
and further with the help of equation (\ref{SDE}), we have the
following derivation

\begin{eqnarray} \label{v_t1}
 dv(t,X_t)&=&\left\{\frac{\partial}{\partial
 t}v(t,X_t)+\frac{1}{2}[Tr(\sigma\sigma^\ast)\nabla^2v](t,X_t)\right.\nonumber
 \\
 && \, \left.+\langle b,\nabla v\rangle(t,X_t)\right\}dt+
 \langle(\sigma^\ast\nabla v)(t,X_t),dB_t\rangle
\end{eqnarray}
since $$\langle\nabla v(t,X_t),\sigma(t,X_t)dB_t\rangle
=\langle\sigma^\ast(t,X_t)\nabla v(t,X_t),dB_t\rangle\, .$$ Now
comparing (\ref{v_t2}) and (\ref{v_t1}) and using the uniqueness of
Doob-Meyer's decomposition of continuous semi-martingale, we
conclude that the coefficients of $dt$ and $dB_t$ must coincide,
respectively, namely
$$(\sigma^{-1}b)(t,X_t)=(\sigma^\ast\nabla v)(t,X_t)$$
and
$$\frac{1}{2}|(\sigma^{-1}b)(t,X_t)|^2=
\frac{\partial}{\partial
t}v(t,X_t)+\frac{1}{2}[Tr(\sigma\sigma^\ast\nabla^2
v)](t,X_t)+\langle b,\nabla v\rangle(t,X_t)$$ holds for all $t>0$.
Since the SDE (\ref{SDE}) is non-degenerate, the support of
$X_t,t\in[0,\infty)$ is the whole space $\R^d$. Hence, the following
two equalities
\begin{equation}\label{eqZ_t1} (\sigma^{-1}b)(t,x)=(\sigma^\ast\nabla)
v(t,x)\end{equation} and
\begin{equation}\label{eqZ_t}
\frac{1}{2}|(\sigma^{-1}b)(t,x)|^2=\frac{\partial}{\partial
t}v(t,x)+\langle b,\nabla
v\rangle(t,x)+\frac{1}{2}[Tr(\sigma\sigma^\ast\nabla v)](t,x)
\end{equation}
hold on $[0,\infty)\times\R^d$. It is clear that equality
(\ref{eqZ_t1}) is nothing but equality (\ref{eqv}), while by
(\ref{eqv}) the equality (\ref{eqZ_t}) reduces to equation
(\ref{KPZ-type}).

\medskip

Now let us turn to the {\it sufficiency}. We assume that there
exists a $C^{1,2}$ scalar function $v:[0,\infty)\times\R^d\to\R$
solving equation (\ref{KPZ-type}). We specify the drift $b$ of SDE
(\ref{SDE}) via (\ref{eqv}), namely
$$b(t,x)=(\sigma\sigma^\ast\nabla v)(t,x),\quad
(t,x)\in[0,\infty)\times\R^d\, .$$ Combining equality (\ref{eqv})
and equation (\ref{KPZ-type}) with equality (\ref{v_t1}), we have
\begin{eqnarray*}
dv(t,X_t) &=& \big[-\frac{1}{2}|\sigma^\ast\nabla v|^2(t,X_t)
+\langle b,\nabla v\rangle(t,X_t)\big]dt \\
&& \qquad +\langle(\sigma^\ast\nabla v)(t,X_t),dB_t\rangle \\
&=&\frac{1}{2}|\sigma^{-1}b|^2(t,X_t)dt
+\langle(\sigma^{-1}b)(t,X_t),dB_t\rangle\, .\end{eqnarray*} This
clearly implies equality (\ref{RND0}) by taking stochastic
integration. We are done. {\it Q.E.D.}

\subsection{The special case of $d=1$}
In this subsection, we would like to discuss our Theorem \ref{TH1}
on $\R$ --- the simplest case. We start with SDE in one dimension:
\begin{equation}\label{SDE-1D}
dX_t=b(t,X_t)dt+\sigma(t,X_t)dB_t, \quad t\ge0\end{equation} with
the diffusion coefficient satisfies that $\sigma(t,x)\neq0$ for all
$(t,x)\in[0,\infty)\times\R$. In this case we have
$$\gamma(t,x)=-\frac{b(t,x)}{\sigma(t,x)}.$$
We set
\begin{equation}\label{ratio-1D}
u(t,x):=\frac{b(t,x)}{\sigma^2(t,x)}=-\frac{\gamma(t,x)}{\sigma(t,x)},\quad
(t,x)\in[0,\infty)\times\R.\end{equation} With the assumption on
$\gamma$ for the Girsanov theorem, we can rephrase our Theorem
\ref{TH1} in a slightly more concise manner

\begin{Theorem} \label{TH1-1D}
Let $v:[0,\infty)\times\R\to\R$ be $C^1$ with respect to the first
variable and $C^2$ with respect to the second variable. Then
\begin{equation}\label{RND-1D}
v(t,X_t)=\frac{dQ_t}{dP}=\exp\left(v(0,X_0)-\int^t_0\frac{b(s,X_s)}{\sigma(s,X_s)}dB_s
-\frac{1}{2}\int^t_0\big|\frac{b(s,X_s)}{\sigma(s,X_s)}\big|^2ds\right)
\end{equation} if and only if there exists a $ C^1$-function $\Phi:
\R\to\R$ such that
$$b(t,X_t)=\Phi(u(t,X_t)),\quad \forall t \geq 0$$
and the function $u(t,x)$ satisfies the following generalized
Burgers equation (again time-reversed)
\begin{equation} \label{struc equ}
    \frac{\partial}{\partial t}u (t,x)=-\frac{1}{2}\frac{\partial^2}{\partial x^2 }
 \Psi_1(u(t,x))-\frac{1}{2}\frac{\partial}{\partial x} \Psi_2(u(t,x))\\
\end{equation}
where
\[
    \Psi_1(r):=\int \frac{\Phi (r)}{r}dr, \quad \Psi_2(r):=r \Phi (r), \quad r\in\R\, .
\]
\end{Theorem}

\medskip
\par\noindent
{\it Proof} The proof is in the same manner as the proof to Theorem
\ref{TH1} together with the combination of the introduction and
properties of the functions $\Psi_1$ and $\Psi_2$. We omit the whole
derivation here. {\it Q.E.D.}

\begin{Remark}
We would like to point out that the function $u$ defined in formula
(\ref{ratio-1D}) has the following explanation. Actually, from It\^o
formula, one may see that the square of Brownian motion has certain
contribution to the drift of the stochastic differential equation.
So the composition $u(t,X_t)$ of the function $u$ with the process
$X_t$ may characterize the proportion of the drift part with respect
to the diffusion part in equation (\ref{SDE-1D}). Surprisingly, this
function $u$ satisfies the nonlinear parabolic PDE of Burgers type
(\ref{struc equ}).
\end{Remark}

Our PDE (\ref{struc equ}) covers much more classes of specific
nonlinear PDEs. Now let us give several examples to explicate this
point.

\begin{Example}
Give a constant $\sigma>0$. Let $b(t,x)=\sigma^2 u (t,x)$ and
$\sigma(t,x)\equiv\sigma$, our SDE (\ref{SDE}) then becomes
\[
    dX_t=\sigma^2 u (t,X_t)dt+ \sigma dB_t.
\]
The $C^1$-function $\Phi$ is simply given by $\Phi(r)=\sigma^2r$ and
the corresponding PDE (\ref{struc equ}) is a classical Burgers
equation (time-reversed)
$$\frac{\partial}{\partial t}u
(t,x)=-\frac{\sigma^2}{2}\frac{\partial^2}{\partial x^2 }
u(t,x)-\sigma^2 u(t,x)\frac{\partial}{\partial x} u(t,x).
$$
This example recovers the main result obtained in Hodges and
Carverhill \cite{HC1993}. Moreover, our Theorem \ref{TH1-1D} also
covers the results obtained in Hodges and Liao \cite{HL2004} and in
Stein and Stein \cite{Stein}.

\end{Example}

The next example shows that our PDE (\ref{struc equ}) can be a
porous media type partial differential equation.
\begin{Example}
We fix $m\in\N$. Let $a(t,x)= m[u(t,x)]^m$ and $b(t,x)=
\sqrt{m}[u(t,x)]^{\frac{m-1}{2}}$, our SDE (\ref{SDE-1D}) then
becomes
$$dX_t=m [u(t,X_t)]^mdt+\sqrt{m}[u(t,X_t)]^{\frac{m-1}{2}}dB_t.
$$
The $C^1$-function $\Phi$ is then given by $\Phi(r)=mr^m$ and the
corresponding PDE (\ref{struc equ}) is a porous media type nonlinear
PDE
\begin{equation*}
\frac{\partial}{\partial
t}u(x,t)=-\frac{1}{2}\frac{\partial^2}{\partial x^2 } u^m(t,x)-m
\frac{\partial}{\partial x} u^{m+1}(t,x).
\end{equation*}
\end{Example}

Our third example is to show that in the time-homogeneous case in
the sense that $b$ and $\sigma$ are functions of the variable
$x\in\R$ only, the corresponding PDE (\ref{struc equ}) then
determines a harmonic function.
\begin{Example}
Let $b(t,x)=b(x)$ and $\sigma(t,x)=\sigma(x)$, our SDE
(\ref{SDE-1D}) then reads as follows
$$dX_t=b(X_t)dt+ \sigma(X_t)dB_t$$
and the corresponding PDE (\ref{struc equ}) is a second order
elliptic equation for harmonic functions
\begin{equation*}
     \frac{\partial^2}{\partial x^2}\Psi_1(u(x))+\frac{\partial}{\partial x}\Psi_2(u(x))=0\\
\end{equation*}
where
\[
    \Psi_1(r)=\int \frac{\Phi (r)}{r}dr, \quad \Psi_2(r)=r \Phi (r),\quad r\in\R.
\]
\end{Example}

\section{Extension to differential manifolds}
In this final section, we extend our Theorem \ref{TH1} for SDEs on a
general connected complete differential manifold. We start with the
following observation. In the situation of SDE (\ref{SDE}) on
$\R^d$, let
$g_t=(g^{ij}_t(\cdot)):=(\sigma\sigma^\ast)^{-1}(t,\cdot)$. Then we
have a time-dependent metric on $\R^d$ defined as follow
$$\langle x,y\rangle_{g_t}:=\sum^d_{i,j=1}g^{ij}_tx_iy_j=\langle g_tx,y\rangle,\quad x,y\in\R^d.$$
Let $\nabla_{g_t}$ and $\Delta_{g_t}$ be the associated gradient and
Laplacian, respectively. Then the generator for the solution to SDE
(\ref{SDE}) can be reformulated as follows (cf. \cite{IkWa})
$$L_tf=\frac{1}{2}\Delta_{g_t}f+\langle\tilde{b}(t,\cdot),\nabla_{g_t}f\rangle_{g_t}$$
for some smooth function $\tilde{b}:[0,\infty)\times\R^d\to\R^d$.
From this point of view, we intend to extend our Theorem \ref{TH1}
to a general connected complete differential manifold.

Now let $M$ be a $d$-dimensional connected complete differential
manifold with a family of Riemannian metrics
$\{g_t\}_{t\in[0,\infty)}$, which is smooth in $t\in[0,\infty)$.
Clearly $(M,g_t)$ is a Riemannian manifold for each
$t\in[0,\infty)$. Let $\{b(t,\cdot)\}_{t\in[0,\infty)}$ be a family
of smooth vector fields on $M$ which is smooth in $t$ as well. Let
$\nabla_{g_t}$ and $\Delta_{g_t}$ denote the gradient and Laplacian
operators induced by the metric $g_t$, respectively. Then the
diffusion process on $M$ generated by the operator
$$L_t:=\frac{1}{2}\Delta_{g_t}+b(t,\cdot)$$
can be constructed by solving the following SDE on $M$
\begin{equation}\label{SDE_M}
dX_t=b(t,X_t)dt+\Phi_t\circ dB_t\end{equation} where
$\{B_t\}_{t\in[0,\infty)}$ is the $d$-dimensional Brownian motion,
$\circ d$ stands for the Stratonovich differential, and $\Phi_t$ is
the horizontal lift of $X_t$ onto the frame bundle $O_t(M)$ of the
Riemannian manifold $(M,g_t)$, namely, $\Phi_t$ solves the following
equation
$$d\Phi_t=H_{\Phi_t}\circ dX_t,$$
with $H_{\Phi_t}:T(M)\to O_t(M)$ being the horizontal lift. Here
$T_t(M)$ denotes the tangent bundle of $M$.

The following result is an extension of our Theorem \ref{TH1} to
$M$.

\begin{Theorem} \label{TH2}
Let $v:[0,\infty)\times M\to\R$ be $C^1$ with respect to the first
variable and $C^2$ with respect to the second variable. Then

\begin{equation}\label{RND-M}
v(t,X_t)=v(0,X_0)+\frac{1}{2}\int^t_0\big|b(s,X_s)\big|^2_{g_t}ds
+\int^t_0\langle(\Phi^{-1}_sb(s,X_s),\circ dB_s\rangle_{g_t}
\end{equation}
holds if and only if
\begin{equation}\label{eqv-M}
b(t,x)=(\nabla_{g_t}v)(t,x),\quad (t,x)\in[0,\infty)\times
M\end{equation} and the following time-reversed KPZ type equation
\begin{equation}\label{KPZ-type-M}
\frac{\partial}{\partial
t}v(t,x)=-\frac{1}{2}\left[(\Delta_{g_t}v)(t,x) +|\nabla_{g_t}
v|^2_{g_t}(t,x)\right]
\end{equation}
hold, where $|z|^2_{g_t}:=\langle z,z\rangle_{g_t}$ for any vector
$z$ on $M$.
\end{Theorem}

\medskip
\par\noindent
{\it Proof} By (\ref{RND-M}), we have
\begin{equation}\label{dv-M}
dv(t,X_t)=\frac{1}{2}|b(t,X_t)|^2_{g_t}dt+\langle\Phi^{-1}_tb(t,X_t),\circ
dB_t\rangle_{g_t}\, .\end{equation} On the other hand, by
(\ref{SDE_M}) and the It\^o formula, we get
\begin{equation}\label{dv-MIto}
dv(t,X_t)=\langle\Phi^{-1}_t\nabla_{g_t}v(t,X_t),\circ
dB_t\rangle_{g_t}+\left\{\frac{1}{2}\Delta_{g_t}v+\langle
b,\nabla_{g_t}v\rangle_{g_t}\right\}(t,X_t)dt\, .
\end{equation}
Now combining (\ref{dv-MIto}) with (\ref{dv-M}), we arrive the
following
$$\nabla_{g_t}v(t,X_t)=b(t,X_t)$$
and
$$\left\{\frac{1}{2}\Delta_{g_t}v+\langle
b,\nabla_{g_t}v\rangle_{g_t}\right\}(t,X_t)=\frac{1}{2}|b(t,X_t)|^2_{g_t}\,
.$$ Since $\{X_t\}_{t\in[0,\infty)}$ is supported by the whole
manifold, the above two equalities imply (\ref{eqv-M}) and
(\ref{KPZ-type-M}), respectively.

On the other hand, combining (\ref{eqv-M}) and (\ref{KPZ-type-M})
with (\ref{dv-MIto}), we obtain (\ref{dv-M}), which implies
(\ref{RND-M}) by stochastic integration. This completes the proof.
{\it Q.E.D.}

\bigskip
\par\noindent
\textbf{Acknowledgment}: The second named author is supported by
WIMCS and NNSFC (10721091). The third named author would like to
thank Beijing Normal University for warm hospitality and stimulating
working atmosphere. The forth named author is grateful to the
support of a departmental Ph D scholarship at Department of
Mathematics, Swansea University.

\end{document}